\documentclass[preprint,3p,a4paper]{elsarticle}
\usepackage{mathrsfs}
\usepackage{amsfonts,amsmath,amssymb,amsthm}

\usepackage{latexsym}
\usepackage{tikz}
\usepackage{hyperref}
\usepackage{geometry}
 \theoremstyle{plain}

  \newtheorem{thm}{Theorem}[section]
  \newtheorem{cor}[thm]{Corollary}
  \newtheorem{lem}[thm]{Lemma}
  \newtheorem{prop}[thm]{Proposition}
  \newtheorem{exmp}[thm]{Example}
  \newtheorem{rem}[thm]{Remark}
  
  \newtheorem{defn}[thm]{Definition}

\newcommand{\da}{\downarrow \hspace{-2pt}}

\newcommand{\ua}{\uparrow \hspace{-2pt}}

\newcommand{\UUa}[0]{\makebox[1ex][l]{\lower.15ex
                                 \hbox{$\uparrow$}}\kern-1ex\lower-.15ex
                                 \hbox{$\uparrow$}}

\newcommand{\DDa}[0]{\makebox[1ex][l]{\lower.15ex
                                 \hbox{$\downarrow$}}\kern-1ex\lower-.15ex
                                 \hbox{$\downarrow$}}

\begin{document}
\begin{frontmatter}
\title{\bf   $\mathcal{S}^*$-convergence and locally hypercompact spaces
\tnoteref{t1}}
\tnotetext[t1]{Research supported by NSF of China (No. 11871353).}

\author{Yuxu Chen}
\ead{yuxuchen@stu.scu.edu.cn}
\author{Hui Kou\corref{cor}}
\ead{kouhui@scu.edu.cn}

\cortext[cor]{Corresponding author}
\address{College of Mathematics, Sichuan University, Chengdu 610064, China}

\begin{abstract}

We introduce the notion of $\mathcal{S}^*_X$-convergence on a $T_0$ topological space $X$ as generalization of $S$-convergence in domain theory, and define finitely approximated spaces. Monotone determined spaces are natural topological extensions of dcpos. The main results are: (1) For a monotone determined space $X$,  $\mathcal{S}^*_X$-convergence is topological iff $X$ is a locally hypercompact space.  
(2) For a $T_0$ space $X$, $\mathcal{S}^*_X$-convergence is topological iff $X$ is a finitely approximating space. (3) If the Lawson topology on a monotone determined space $X$ is compact, then $X$ is a dcpo endowed with the Scott topology.
(4) A characterization of the closure of a subset of a locally hypercompact space is given.


  \vskip 2mm

\end{abstract}

\end{frontmatter}

\section{Introduction}

With the development of domain theory, it is found to be strongly connected with general topology \cite{CON03,GAU}.  The notion of $\mathcal{S}$-convergence \cite{CON03} was introduced to characterize continuous domains, which reveals the deep connection between order and topology.

\begin{thm}\rm \cite[Theorem \uppercase\expandafter{\romannumeral 2}-1.9]{CON03}
For a dcpo $P$, $P$ is a continuous domain iff the $\mathcal{S}$-convergence on $P$ is topological.
\end{thm}
At the 6th International Symposium in Domain Theory, J.D. Lawson called to develop the core of domain theory directly in $T_0$ spaces instead of posets and showed that some results in
domain theory can be lifted from the context of posets to $T_0$ spaces.
 Recently, Zhao and Lu extended the lim-inf convergence to $T_0$ enriched closure spaces in \cite{ZHAO2022}. Zhang, Bao and Xu \cite{Z2022} defined the notion of $\mathcal{S}$-convergence on $T_0$ spaces, and give out the topological generalization of the above theorem (see \cite[Theorem 4.17]{Z2022}).





In \cite{Li2013}, the notion of $\mathcal{S}^*$-convergence was introduced to characterize the notion of quasicontinuous domains.
Monotone determined spaces were introduce in \cite{ERNE2009} and were shown to be very appropriate topological extensions of dcpos. The notion of monotone determined spaces is equivalent to directed spaces in \cite{ERNE2009}. Many classical structures such as c-spaces, locally hypercompact spaces are all monotone determined spaces. Moreover, by introducing the notion of a generalized way-below relation on a monotone determined space, c-spaces can be viewed as continuous monotone determined spaces and locally hypercompacts spaces can be viewed as quasicontinuous monotone determined spaces \cite{FK2017,XK2022,ZHA2021}, their roles in monotone determined spaces just like continuous dcpos and quasicontinuous dcpos in dcpos. In this paper, we define the notion of   $\mathcal{S}^*_X$-convergence on a topological space $X$. We will show that for a monotone determined space $X$, $\mathcal{S}^*_X$-convergence is topological iff $X$ is a locally hypercompact space (quasicontinuous space). Moreover, We will introduce the notion of finitely approximated spaces and show that a $T_0$ space $X$ is a finitely approximated space iff $X$ endowed with the monotone determined topology is a locally hypercompact space, iff $\mathcal{S}^*_X$-convergence is topological.

Lawson topology is important in domain theory. The liminf convergence on a continuous dcpo is topological \cite{CON03} and the topology generated by liminf convergence is equal to the Lawson topology. In \cite{Xi2016}, it was shown that a dcpo endowed with the Scott topology is Lawson compact iff it is well-filtered, compact and coherent. We extend the notion of Lawson topology and the liminf topology to a monotone determined space. We show that for a locally hypercompact space $X$, the quasi liminf topology is equal to the Lawson topology. Besides, we show that for a monotone determined space, if the Lawson topology is compact, then $X$ is a dcpo endowed with the Scott topology. Finally, for locally hypercompact spaces, we characterize the closure of any subset of $X$ by the limits of all its directed subsets.

\section{Preliminlar}
Let $P$ be a poset. For $A\subseteq P$, we set $\da A=\{x\in P: \exists a\in A,  \ x\leq a\}$,\ $\ua A =\{x\in P: \exists a\in P, \ a\leq x\}$.\ $A$ is called a lower or upper set, if $A=\ \da A$ or $A=\ \ua A$ respectively. For an element $a\in P$, we use $\da a$ or $\ua a$ instead of $\da\{a\}$ or $\ua \{a\}$, respectively.  

For $x,y\in P$, we say that $x$ is way-below $y$, denoted by $x\ll y$, if for any directed subset $D$ of $P$ with the supremum $\bigvee D$ existing, $y\leq \bigvee D$ implies $x\leq d$ for some $d\in D$. $P$ is called a continuous poset if $\{a\in P: a\ll x\}$ is directed and has $x$ as its supremum for all $x\in P$. 

Topological spaces will always be supposed to be $T_0$. For a topological space $X$, its topology is denoted by $\mathcal{O}(X)$ or $\tau$. The partial order $\sqsubseteq$ defined on $X$ by $x\sqsubseteq y \Leftrightarrow  x\in \overline{\{y\}}$
is called the  specialization order, where $\overline{\{y\}}$ is the closure of $\{y\}$.\ From now on, all order-theoretical statements about  $T_0$ spaces, such as upper sets, lower sets, directed sets, and so on,  always refer to  the specialization order $\sqsubseteq $.\ For any net $\xi = (x_i)_{i \in I}$,\ where $I$ is a directed set,\ we write it $(x_i)_I$ for short.\ Given any $x\in X$, $(x_i)_I$ is called converging to $x$,\ denoted by $(x_i)_I  \rightarrow x $,\ if $\{ x_i \}_I$ is eventually in every open neighborhood of $x$.\ The interior of $A$ is denoted by $A^\circ$ or $int (A)$. The closure of $A$ is denoted by $\overline{A}$ or $cl(A)$.

\vskip 3mm
For any set $X$, we denote $\Phi(X)$ to be the class of all nets in $X$. A convergence class  $\mathcal{E}$ in $X$ is a relation between $\Phi X$ and $X$, i.e., $\mathcal{E}$ is a subclass of $\{ (\xi,x): \xi \text{ is a net in }X,\ x \in X\}$. A topological space determined by $\mathcal{E}$,
denoted by $(X,\mathcal{E}(X))$ or $\mathcal{E}X$ for short, can be defined as follows.
 $$U \in \mathcal{E}(X)\ \Longleftrightarrow\   \forall  (\xi,x) \in \mathcal{E},\  x \in U \text{ implies that } \xi \text{ is eventually in } U. $$
On the other hand, given any topological space $ (X,\tau) $,\ we define a convergence class $\mathcal{C}_{(X,\tau)} $\  to be $\{(\xi,x): \xi \to x \text{ in } X\}$.\ 
Obviously,\ for any topological space $(X,\mathcal{E}(X))$,\ we have $\mathcal{E} \subseteq \mathcal{C}_{(X,\mathcal{E}(X))}$, i.e., $\xi \to x$ in $\mathcal{E}X$ for any $(\xi,x) \in \mathcal{E}$.

We say that a convergence class $\mathcal{E}$ on a set $X$ is topological if there is a topology $\tau$ on $X$ such that $\mathcal{E} = \mathcal{C}_{(X,\tau)}$. It is easy to see that $\mathcal{E}$ is topological iff $\mathcal{E}  = \mathcal{C}_{\mathcal{E}X}$. For any two convergence class $\mathcal{S}$ and $\mathcal{E}$, if $\mathcal{S} \subseteq \mathcal{E}$, then  $\mathcal{E}(X) 
\subseteq \mathcal{S}(X)$. We also call a convergence class $\mathcal{S}$ as $\mathcal{S}$-convergence. 

\vskip 3mm

Let $(X,\mathcal{O}(X))$ be a $T_0$ space. Every directed subset $D\subseteq X$ can be  regarded as a monotone net $(d)_{d\in D}$. Set
$DS(X)=\{D\subseteq X: D \ {\rm is \ directed}\}$ to be the family of all directed subsets of $X$. For an $x\in X$, we use  $D\rightarrow x$  to denote  that $x$ is a limit of $D$, i.e., $D$ converges to $x$ with respect to $\mathcal{O}(X)$. Then the following result is obvious.

\begin{lem}  Let $X$ be a $T_0$ space. For any $(D,x)\in DS(X)\times X$, $D\rightarrow x$ if and only if $D\cap U\not=\emptyset$ for any  open neighborhood of $x$.
\end{lem}

Set $DLim(X)=\{(D,x)\in DS(X)\times X:  \   D\rightarrow x \}$ to be the set of all pairs of directed subsets and their limits in  $X$. Then $(\{y\},x)\in DLim(X) $ iff $x\sqsubseteq y$ for all $x,y\in X$. 

Let $X$ be a $T_0$ space. A subset $U\subseteq X$ is called monotone determined open  if for all $(D,x)\in DLim(X)$, $x\in U$ implies $D\cap U\not=\emptyset$. Obviously, every open set of $X$ is monotone determined open. Set
$\mathcal{D}(X) =\{U\subseteq X: U \ {\rm is} \  {\rm monotone\ determined\ } {\rm open}\},$
then $\mathcal{O}(X)\subseteq \mathcal{D}(X)$.

\vskip 3mm

\begin{defn}\rm \cite{L2016,YK2015} A  topological space $X$ is said to be a monotone determined space if it is $T_0$ and  every monotone determined open set is open; equivalently, $ \mathcal{D}(X)=\mathcal{O}(X)$.
\end{defn}

The following are some basic properties of monotone determined spaces. Given any space $X$, we denote $\mathcal{D}X$ to be the topological space $(X,\mathcal{D}(X))$.
\vskip 3mm

\begin{thm}\rm \cite{L2016,YK2015} \label{convergence}
 Let $X$ be a $T_0$ topological space.
\begin{enumerate}
\item[(1)] For all $U\in \mathcal{D}(X)$, $U=\ \ua U$.
\item[(2)] $X$ equipped with topology $\mathcal{D}(X)$ is a $T_0$ topological space such that  $\sqsubseteq_d =\sqsubseteq$, where $\sqsubseteq_d$ is the specialization order relative to $\mathcal{D}(X)$.
\item[(3)] For a directed subset $D$ of $X$, $D\rightarrow x$ iff $D\rightarrow_d x$ for all $x\in X$, where $D\rightarrow_d x$ means that $D$ converges to $x$ with respect to the topology $\mathcal{D}(X)$.
\item[(4)] $\mathcal{D}X$ is a monotone determined space.
\end{enumerate}
\end{thm}



Every dcpo endowed with the Scott topology is a monotone determined space. Every poset endowed with the weak Scott topology \cite{ERNE2009} is a monotone determined space. The notion of a monotone determined space is a very natural topological extensions of dcpos, and the category of all monotone determined spaces with continuous maps as morphisms is cartesian closed \cite{FK2017,L2016,XK2022,YK2015}.

\vskip 2mm

A space $X$ is called a c-space if $\forall x \in U \in \mathcal{O}(X), \exists y \in X \text{ s.t. } x \in (\ua y)^\circ \subseteq\ \ua y \subseteq U$. A space $X$ is called a locally hypercompact space if $\forall x \in U \in \mathcal{O}(X), \exists F \subseteq_f X \text{ s.t. } x \in (\ua F)^\circ \subseteq\ \ua F \subseteq U$, where $F \subseteq_f X$ means that $F$ is a finite subset of $X$. C-spaces and locally hypercompact spaces can be characterized as continunous and quasicontinuous monotone determined spaces respectively. 

Given any set $X$, we use $F \subseteq_{f} X$ to denote that $F$ is a finite subset of $X$.\ Given any two subsets $G,H \subseteq X$,\ we define $G \leq H$ iff $\ua H \subseteq \ua G$.\ A family of finite sets is said to be directed if given any two $F_{1}, F_{2}$ in the family,\ there exists an $F$ in the family such that $F_{1},F_{2} \leq F$.\

\begin{defn}\rm \cite{FK2017,ZHA2021}   \label{d-def}
Let $X$ be a monotone determined space, $y \in X$ and $G,H \subseteq X$.\ We say that $G$ $d$-approximates $H$, denoted by $G \ll_{d} H$, if for every directed subset $D\subseteq X$, $D \to y$ for some $y\in H$ implies $D\ \cap \ua\! G \neq \emptyset$. We write $G \ll_{d} x$ for $G \ll_{d} \{ x \}$. 
\end{defn}

\vskip 3mm

Let $X$ be a $T_{0}$ topological space and $\mathcal{M}$ be a directed family of finite subsets of $X$.\ Then $F_{\mathcal{M}} = \{A \subseteq X : \exists M \in \mathcal{M}, \ua M \subseteq A\}$ is a filter under the inclusion order.\ We say that $\mathcal{M}$ converges to $x$, denoted by $\mathcal{M} \to x$, if $F_{\mathcal{M}}$ converges to $x$,\ i.e.,\ for any open neighbourhood $U$ of $x$,\ there exists some $M \in \mathcal{M}$ such that $M \subseteq U$.\

\vskip 3mm

\begin{defn} \rm \cite{FK2017,ZHA2021} \label{def-quasi}
A topological space $X$ is called a quasicontinuous space if it is a monotone determined space such that for any $x \in X$,\ the family $fin_{d}(x) = \{F:\ F  $ is finite,\ $F \ll_{d} x \}$ is a directed family and converges to $x$.
\end{defn}

\vskip 3mm

For any space $X$ and any finite subset $F$, denote $\Uparrow_d F$ to be $\{x \in X: F \ll_d x\}$.

\vskip 3mm

\begin{thm}\rm \cite{FK2017,ZHA2021} \label{d-quasicontinuity}
$X$ is a quasicontinuous space iff $X$ is a locally hypercompact space. Let $X$ be a quasicontinuous space. The following statements hold.
\begin{enumerate}
\item[(1)] Given any $H \subseteq X$ and  $y\in X$, $H \ll_d y$ implies $H \ll_d F\ll_d y$ for some finite subset $F \subseteq_{f} X$.

\item[(2)] Given any $F \subseteq_{f} X $, $\Uparrow_d F =(\ua F)^{\circ}$. Moreover, $\{\Uparrow_d F: F \subseteq_{f} X\}$ is a base of the topology of $X$.

\item[(3)] Given any $G,H \subseteq X$, the following are equivalent:
\begin{enumerate}
\item[(i)] $G\ll_d H$;
\item[(ii)] $H \subseteq (\ua G)^{\circ}$.
\item[(iii)] For any net $(x_j)_J\subseteq X$, $(x_j)\rightarrow y, y \in H $ implies $\ua G \cap  (x_{j})_{J} \neq \emptyset$.

\end{enumerate}
\end{enumerate}

\end{thm}

\vskip 3mm

Therefore, the notion of locally hypercompact spaces and quasicontinuous coincides. As special cases of quasicontinuous spaces, the notion of continuous spaces agrees with the notion of c-spaces.

Given any topological space $X$ and $x,y \in X$,
 we say that $x$ $d$-approximates $y$, denoted by $x\ll_d y$, if $\{x\} \ll_d \{y\}$ in Definition \ref{d-def}. For any monotone determined space $X$ and $x\in X$, we denote 
$ \UUa_d x = \{y\in X: x\ll_d y\}$. 
Let $X$ be a  monotone determined space.
$X$ is called a continuous space if $\DDa_dx$ is directed and $\DDa_dx\rightarrow x$ for all $x\in X$.
\vskip 3mm

\begin{thm}\rm \cite{XK2022}
Let $X$ be a monotone determined space. Then $X$ is a continuous space iff $X$ is a c-space iff for any $x \in X$ there exists a directed set $D \subseteq \DDa_d x$ such that $D \to x$.
\end{thm}
\vskip 3mm






Continuous dcpos endowed with the Scott topology and $s_2$-continuous posets endowed with the weak Scott topology ($\sigma_2$ topology) can both be viewed as special continuous spaces \cite{ZHA2021}. 
From now on, we drop the subscript $d$ in these notations when there is no confusion can arise.

In Definition \ref{def-quasi}, if we do not require the space $X$ to be a monotone determined space, then we call it a finitely approximated space. There exists a space that is a finitely approximated space but not monotone determined hence not a quasicontinuous space.

\vskip 3mm
\begin{defn} \rm 
A $T_0$ space is called a finitely approximated space if for any $x \in X$, $fin(x) = \{F: F \subseteq_f X, F \ll_d x\}$ is a directed family and converges to $x$.
\end{defn}
\vskip 3mm

\begin{exmp} \rm  \label{exand}
Let $\mathbb{N}$ be the set of natural numbers. Denote $\mathbb{N}^\top$ the flat domain, i.e., the poset with carrier set $\mathbb{N} \cup \{\top\}$ and $x \leq y$ iff $y = \top$ or $x=y$. 
It is easily seen that $\mathbb{N}^\top$ is a dcpo. 
Considering the upper topology on $\mathbb{N}^\top$, it is easy to verify the following: 
\begin{enumerate}
\item[(1)]  $(\mathbb{N}^\top, v(\mathbb{N}^\top))$ is a locally compact sober space;
\item[(2)]  $(\mathbb{N}^\top, v(\mathbb{N}^\top))$ is a finitely approximated space.
\item[(3)]   $(\mathbb{N}^\top, v(\mathbb{N}^\top))$ is not a monotone determined space.
\end{enumerate}
\end{exmp}

\vskip 3mm

The following  properties of finitely approximated spaces are routine to verify in a same way as quasicontinuous domains \cite{CON03} and quasicontinuous spaces \cite{FK2017,ZHA2021}, we omit the proof. 

\begin{prop}\rm  \label{interpolation}
Let $X$ be a finitely approximated space. 
\begin{enumerate}
\item[(1)] $F \ll_d G$ implies $H \leq G$.
\item[(2)] $F_1 \leq F \ll_d G \leq G_1$ implies $F_1 \ll_d G_1$.

\item[(3)] Given any $H \subseteq_f X$ and $y \in X$ such that $H \ll_d y$, then there exists some finite subset $F \subseteq_f X$ such that $H \ll_d F \ll_d y$.
\end{enumerate}
\end{prop}

\section{Characterizing continuity of spaces by $\mathcal{S}_X$-convergence}

 As mentioned in the introduction, $\mathcal{S}$-convergence on a dcpo is topological if and only if the dcpo is a continuous domain \cite{CON03}.\ The notion of $\mathcal{S}^*$-convergence on a dcpo was introduced in \cite{Li2013} and it was shown that $\mathcal{S}^*$-convergence in a dcpo is topological if and only if the dcpo is a quasicontinuous domain. Later in \cite{R2016}, the notion of $\mathcal{S}^*$-convergence were modified to characterize  $s_2$-quasicontinuous posets.

In the previous section, we see that both dcpos endowed with the Scott topology and posets endowed with the weak Scott topology are monotone determined spaces. Quasicontinuous domains endowed with the Scott topology and $s_2$-quasicontinuous posets endowed with the weak Scott topology are both special quasicontinuous spaces. 
We will define the notion of  $\mathcal{S}^*_X$-convergence on a topological spaces $X$. We show that $\mathcal{S}^*_X$-convergence on a monotone determined space is topological iff $X$ is and a quasicontinuous space. Moreover, $\mathcal{S}^*_X$-convergence on a $T_0$ spaces is topological iff $\mathcal{D}X$ is a quasicontinuous space iff $X$ is a finitely approximated space.


\begin{defn} \rm  
Let $X$ be a $T_0$ topological space. Given any net $(x_i)_{i \in I}$ of $X$, we say that $y$ is an eventual lower bound of $(x_i)_{i \in I}$ if there exists a $k \in I$ such that $ y \leq x_i$ for any $i \geq k $. Denote $EL((x_i)_{i \in I})$ the set of eventual lower bounds of $(x_i)_{i \in I}$.

\end{defn}
\vskip 3mm

\begin{defn} \rm 

We define a convergence class $\mathcal{S}_X$ on $X$ as follows: $((x_i)_{i\in I},x) \in \mathcal{S}_X$ iff there exists a directed subset $D \subseteq EL((x_i)_{i\in I})$ such that $D \to x$ with respect to $\mathcal{O}(X)$.
\end{defn}

By definition, $\mathcal{S}_X$-convergence is topological iff $\mathcal{S}_X$-convergence contains all pairs of convergent nets and limits in $\mathcal{S}(X)$, where $\mathcal{S}(X)$ denotes the topology on $X$ determined by $\mathcal{S}_X$, i.e., $U \in \mathcal{S}(X)$ iff for all $ ((x_i)_{i \in I},x) \in \mathcal{S}_X$ and  $x \in U$, $(x_i)_{i \in I}$ is eventually in $U$.

\vskip 3mm

In \cite[Proposition 3.3]{Z2022}, the following statement is in fact proved.

\begin{lem}\rm \cite{Z2022} \label{s-d-topology}
Let $X$ be a $T_0$ space. Then $U \in \mathcal{S}(X)$ iff $U \in \mathcal{D}(X)$.

\end{lem}



\vskip 3mm

\begin{prop}\rm \cite{Z2022}  \label{continuous topological}
Let $X$ be a continuous space. Then $\mathcal{S}_X$ is topological.
\end{prop}

\vskip 3mm

When focusing on monotone determined spaces, the following holds. 
\vskip 3mm
\begin{prop}\rm   \label{topological continuous}
Let $X$ be a monotone determined space. If $\mathcal{S}_X$-convergence is topological, then $X$ is a continuous space.

\end{prop} 
\proof By Lemma \ref{s-d-topology}, $\mathcal{S}(X) =\mathcal{D}(X)$. If $\mathcal{S}_X$-convergence is topological, then $((x_i)_{i\in I},x) \in \mathcal{S}_X$ iff  $(x_i)_{i\in I} \to x$ with respect to $ \mathcal{D}(X) = \mathcal{O}(X)$. Let $x \in X$, define
\[I = \{(U,a) \in \mathcal{N}(x) \times L: a \in U\}, \]
where $\mathcal{N}(x)$ consists of all open sets containing $x$. Define an order on $I$ as follows:  $(U,a) \leq (V,b)$ iff $V$ is a subset of $U$. For each $i = (U,a)$, let $x_i = a$. Then $(x_i)_{i \in I}$ converges to $x$ relative to $\mathcal{O}(X)$. Thus, $((x_i)_{i \in I},x) \in \mathcal{S}_X$. It follows that there exists a directed set $D \subseteq EL((x_i)_{i \in I})$ such that $D \to x$ relative to $\mathcal{O}(X)$. Given any $d \in D$, there exists a $k = (U,a ) \in I$ such that $(V,b) = j \geq k$ implies $d \leq b$. Since $\forall b \in U, (U,b) \geq (U,a )$, we have $d \leq u, \forall u \in U$. Thus, $x \in U \subseteq  (\ua d)^\circ$. Then $D$ is a directed subset with $x \in (\ua d)^\circ,\forall d \in D,$ and $D \to x$. Thus, $X$ is a continuous space. $\Box$

\vskip 3mm

Combining Proposition \ref{continuous topological} and Proposition \ref{topological continuous}, we gain the following statement.

\vskip 3mm

\begin{thm}\rm   \label{con con}
Let $X$ be a monotone determined space. $\mathcal{S}_X$-convergence is topological iff $X$ is a continuous space.
\end{thm}
\vskip 3mm

Then, the following is a direct corollary.
\vskip 3mm
\begin{cor}\rm \cite{Z2022}
Let $X$ be a $T_0$ space. $X$ is a continuous space iff $\mathcal{S}_X$ is equal to $\mathcal{C}_X$.
\end{cor}
\proof
We need only to show the necessity. Suppose that $\mathcal{S}_X = \mathcal{C}_X$. Then $\mathcal{S}_X$-convergence is topological and $\mathcal{D}(X) = \mathcal{S}(X) = \mathcal{O}(X)$. Thus, $X$ is a monotone determined space. $X$ is a continuous space by Theorem \ref{con con}.
$\Box$
\vskip 3mm

The following example shows that for a $T_0$ space $X$ with $\mathcal{S}_X$ topological, $X$ may not be a monotone determined space hence not a continuous space.

\vskip 3mm
\begin{exmp}\rm  \label{example}
Let $\mathbb{N}$ be the set of natural numbers, $X = (\mathbb{N},co(\mathbb{N}))$ be $\mathbb{N}$ endowed with the cofinite topology and $Y$ be $\mathbb{N}$ endowed with the discrete topology. Then $X,Y$ are both $T_1$. A directed subset of $X$ must be a singleton set \{n\} for some $n \in \mathbb{N}$. If $\{n\} \to x$ for some $x \in \mathbb{N}$ relative to $X$, then $n = x$. Thus, $((x_i)_I,x) \in \mathcal{S}_X$ iff there exists some $n \in \mathbb{N}$ such that $n$ is an eventual lower bound of $(x_i)_I$ and $n = x$, iff $x$ is an eventual lower bound of $(x_i)_I$, iff there exists some $i_0 \in I$ such that $x = x_i$ for all $i \geq i_0$. Then $\mathcal{S}_X$ is topological since it is equal to $\mathcal{C}_X$. However, $X$ is not a monotone determined space.
\end{exmp}
\vskip 3mm

\begin{lem}\rm  \label{dx}
Let $X$ be a $T_0$ space. Then $\mathcal{S}_X = \mathcal{S}_{\mathcal{D}X}$.
\end{lem}
\proof
By Theorem \ref{convergence}, the specialization orders of $X$ and $\mathcal{D}X$ are the same and for a directed subset $D$, $D \to x$ relative to $\mathcal{D}X$ iff $D \to x$ relative to $X$. For a net $(x_i)_I$, the eventual lower bounds in $X$ and $\mathcal{D}X$ are the same. By definition, $\mathcal{S}_X = \mathcal{S}_{\mathcal{D}X}$.
$\Box$
\vskip 3mm

\begin{thm} \rm \label{con-top}
Let $X$ be a $T_0$ space. $\mathcal{S}_X$ is topological iff $\mathcal{D}X$ is a continuous space.
\end{thm}

\proof
By Lemma \ref{dx}, $\mathcal{S}_X$ is topological iff $\mathcal{S}_{\mathcal{D}X} $ is topological. Since $\mathcal{D}X$ is a monotone determined space, $\mathcal{S}_{\mathcal{D}X} $ is topological iff $\mathcal{D}X$ is a continuous space by Theorem \ref{con con}.
$\Box$

\vskip 3mm

Now we investigate the relationships between $\mathcal{S}^*_X$-convergence and quasicontinuous spaces (locally hypercompact spaces).

\vskip 3mm

\begin{defn}\rm 
Let $X$ be a $T_0$ topological space. Given any net $(x_i)_{i \in I}$, we say that a finite set $F$ is a quasi eventual lower bound of $(x_i)_{i \in I}$ if there exists a $k \in I$ such that $\forall i \geq k,  x_i \in\ \ua F $. Denote $QEL((x_i)_{i \in I})$ the set of quasi eventual lower bound of $(x_i)_{i \in I}$.
\end{defn}

\vskip 3mm

\begin{defn} \rm 
Let $X$ be a topological space. We define a convergence class $\mathcal{S}^*_X$ on $X$ as follows: $((x_i)_{i\in I},x) \in \mathcal{S}^*_X$ iff there exists a directed family $\mathcal{F}$ of quasi eventual lower bounds of $(x_i)_{i\in I}$ such that $\mathcal{F} \to  x$ (Recall that $\mathcal{F} \to  x$ means that for any open subset $U$ containing $x$, there exists some $F \in \mathcal{F}$ such that $F \subseteq U$).
\end{defn}

\vskip 3mm

By definition, $\mathcal{S}^*_X$-convergence is topological iff $\mathcal{S}^*_X$ contains all pairs of convergent nets and limits in $\mathcal{S}^*(X)$, where $\mathcal{S}^*(X)$ is the topology on $X$ determined by $\mathcal{S}^*_X$, i.e., $U \in \mathcal{S}^*(X)$ iff for all $ ((x_i)_{i \in I},x) \in \mathcal{S}^*_X$ and $x \in U$, $(x_i)_{i \in I}$ is eventually in $U$. For a monotone determined space, we can characterize its quasicontinuity by that $\mathcal{S}^*_X$-convergence is topological.

\vskip 3mm

\begin{lem}[Rudin's Lemma]\rm\cite{CON03}\label{rudin}
Let $\mathcal{F}$ be a directed family of nonempty finite subsets of a poset $P$. Then there exists a directed set $D \subseteq \bigcup_{F \in \mathcal{F}} F$ such that $D \cap F \not = \emptyset $ for all $F \in \mathcal{F}$.
\end{lem}

\vskip 3mm

$\mathcal{S}^*_X$-convergence are generalized $\mathcal{S}_X$-convergence. It is easy to see that $\mathcal{S}_X \subseteq \mathcal{S}^*_X$ since each element $x$ can be viewed as a singleton set $\{x\}$, and then each directed subset $D$ can be viewed as a directed family $\{\ua d: d \in D\}$. The following example shows that, generally, $\mathcal{S}_X \not = \mathcal{S}^*_X$. 

\vskip 3mm

\begin{exmp}  \rm 
Let $X = \mathbb{N} \cup \{\top,a\}$. Given any $x,y \in X$, define $x \sqsubseteq y$ iff $y = \top$ or $x,y \in \mathbb{N},x \leq y$. Let $X$ endowed with the Scott topology with respect to the order $\sqsubseteq$. Consider the net $(x_n)_\mathbb{N}$ with $x_{2n} = a , x_{2n+1} = n$. Then each $\{a,n\}$ for $n \in \mathbb{N}$ is a quasi eventual lower bound of $(x_i)_I$. Let $\mathcal{F} = \{\{a,n\}\}_{n \in \mathbb{N}}$. $\mathcal{F}$ is a directed family and $\mathcal{F} \to a$. So $((x_i)_I,x) \in \mathcal{S}^*_X$. It is easy to see that there is no eventual lower bound of $(x_i)_I$. Therefore, $((x_i)_I,x) \not \in \mathcal{S}_X$.
\end{exmp}
In the following statement, we will see that although $\mathcal{S}_X$ and $\mathcal{S}^*_X$ may be different, they will determined a same topology which is equal to $\mathcal{D}(X)$.

\vskip 3mm

\begin{lem}\rm  \label{directed-convergence-family}
Let $X$ be a $T_0$ space and $\mathcal{F}$ be a directed family of finite subsets of $X$. Then $\mathcal{F} \to x$ with respect to $X$ iff $\mathcal{F} \to x$ with respect to $\mathcal{D}X$.
\end{lem}
\proof
Suppose that $\mathcal{F} \to  $ with respect to $\mathcal{D}X$. Since the topology of $\mathcal{D}X$ is finer than that of $X$, then $\mathcal{F} \to  $ with respect to $X$.

Conversely, suppose that  $\mathcal{F} \to x $ with respect to $X$. 
 We show that for each $U \in \mathcal{D}(X)$ with $x \in U$, there exists some $F \in \mathcal{F}$ such that $F \subseteq U$. If not, then $F^\prime = F \backslash U$ is nonempty for any $F \in \mathcal{F}$. Let $\mathcal{F}^\prime = \{F^\prime : F \in \mathcal{F}\}$. Since $U$ is an upper set, $\mathcal{F}^\prime$ is a directed family as well.  By Rudin's Lemma, there exists a directed set $D \subseteq \bigcup_{F^\prime \in \mathcal{F}^\prime}  F^\prime$ such that $D \cap F^\prime \not = \emptyset$. Since $\mathcal{F} \to x$ with respect to $\mathcal{O}(X)$, then $D \to x$ relative to $\mathcal{O}(X)$. It follows that $D \to x$ with respect to $\mathcal{D}(X)$ by Theorem \ref{convergence}. Thus, $D \cap U \not = \emptyset$, a contradiction. Thus, there exists an $F \in \mathcal{F}$ such that $F \subseteq U$. Since $F \in QEL((x_i)_I)$, $(x_i)_I$ is eventually in $U$. $\Box$

\vskip 3mm

\begin{lem}\rm  \label{qs-d-topology}
Let $X$ be a $T_0$ space. Then $\mathcal{S}^*(X) = \mathcal{S}(X) = \mathcal{D}(X)$.

\end{lem}

\proof Since $\mathcal{S}_X \subseteq \mathcal{S}^*_X$, $\mathcal{S}^*(X) \subseteq \mathcal{S}(X)$. By Lemma \ref{s-d-topology}, we need only to show that $\mathcal{D}(X) \subseteq \mathcal{S}^*(X)$.

Suppose that $U \in \mathcal{D}(X)$. Given any $((x_i)_{i \in I},x) \in \mathcal{S}^*_X$ with $x \in U$,  we need only to show that $(x_i)_{i\in I}$ is eventually in $U$. By definition of $\mathcal{S}^*_X$, there exists a directed family  $\mathcal{F} \subseteq QEL((x_i)_{i \in I})$ such that $\mathcal{F} \to x$ with respect to $X$. 
 By Lemma \ref{directed-convergence-family}, we have $\mathcal{F} \to x$ with respect to $\mathcal{D}X$. Therefore, there exists some $F \in \mathcal{F}$ such that $F \subseteq U$ and then $(x_i)_{i \in I}$ is eventually in $\ua F$ and hence eventually in $U$. 
 $\Box$

\vskip 3mm

\begin{prop}\rm   \label{qc to}
Let $X$ be a quasicontinuous space. Then 
$\mathcal{S}^*_X$ is topological.

\end{prop}
\proof 
If $((x_i)_{i\in I},x) \in \mathcal{S}^*_X$, then there exists some directed family $\mathcal{F}$ of $QEL((x_i)_{i\in I})$ such that $\mathcal{F} \to x$. Let $G$ be any finite set with $ x \in {\Uparrow} G = (\ua G)^\circ$. Then, there exists some $F \in \mathcal{F}$ such that $\ua F \subseteq {\Uparrow} G $. Since $F \in QEL((x_i)_I)$, $(x_i)_I$ is eventually in ${\Uparrow} G$. By Theorem \ref{d-quasicontinuity}, all ${\Uparrow} G$ for $G \subseteq_f X$ form a base of the topology of $X$. It follows that $(x_i)_I$ converges to $x$ relative $\mathcal{O}(X)$.

Conversely, suppose $(x_i)_{i\in I} \to x$. For each $F \in fin(x) = \{F \subseteq_f X : F\ll x\}$, $(x_i)_I$ is eventually in ${\Uparrow} F = (\ua F)^\circ$. Thus, $F$ is a quasi eventual lower of  $(x_i)_{i \in I}$. Since $fin(x)$ is a directed family and $fin(x) \to x$, $((x_i)_{i\in I} , x) \in \mathcal{S}^*$. $\Box$

\vskip 3mm

By Proposition \ref{continuous topological} and Proposition \ref{qc to}, when $X$ is a continuous space, $\mathcal{S}_X$ and $\mathcal{S}^*_S$ are both topological and equal to $\mathcal{C}_X$. Thus, $\mathcal{S}_X = \mathcal{S}^*_X$.
\vskip 3mm
\begin{cor}\rm 
Let $X$ be a continuous space. Then $\mathcal{S}^*_X = \mathcal{S}_X$.
\end{cor}

\vskip 3mm

\begin{lem}\rm  \label{qua-topo}
Let $X$ be a monotone determined space. If $\mathcal{S}^*_X$ is topological, then $X$ is quasicontinuous.
\end{lem}

\proof
Suppose that $\mathcal{S}^*_X$ is topological. Then $((x_i)_{i\in I},x) \in \mathcal{S}^*_X$ iff  $(x_i)_{i\in I} \to x$ with respect to $ \mathcal{D}(X) = \mathcal{O}(X)$. Let $x \in X$, define
\[I = \{(U,a) \in \mathcal{N}(x) \times L: a \in U\}, \]
where $\mathcal{N}(x)$ consists of all open sets containing $x$. Define an order on $I$ as follows: $(U,a) \leq (V,b)$ iff $V \subseteq U$. For each $i = (U,a)$, let $x_i = a$. Then $(x_i)_{i \in I}$ converges to $x$ relative to $\mathcal{O}(X)$. Thus, $((x_i)_{i \in I},x) \in \mathcal{S}^*_X$. It follows that there exists a directed family $\mathcal{F} \subseteq QEL((x_i)_{i \in I})$ such that $\mathcal{F} \to x$ relative to $\mathcal{O}(X)$. Given any $F \in \mathcal{F}$, there exists a $k = (U,a ) \in I$ such that $(V,b) = j \geq k$ implies $ b \in \ua F $. Since $\forall b \in U, (U,b) \geq (U,a )$, we have $U \subseteq \ua F$. Thus, $x \in (\ua F)^\circ$. Then $\mathcal{F}$ is a directed family of finite sets with $\forall F \in \mathcal{F}, x \in (\ua F)^\circ$ and $\mathcal{F} \to x$. Thus, $X$ is a quasicontinuous space. $\Box$

\vskip 3mm

By Proposition \ref{qc to} and Lemma \ref{qua-topo}, we get the following statement. 

\vskip 3mm

\begin{thm} \rm 
Let $X$ be a monotone determined space. $\mathcal{S}^*_X$-convergence is topological iff $X$ is quasicontinuous.
\end{thm}
\vskip 3mm

\begin{thm}\rm    \label{T-qua}
Let $X$ be a $T_0$ space. $\mathcal{S}^*_X$-convergence is topological iff $\mathcal{D}X$ is a quasicontinuous space.
\end{thm}
\proof
By Lemma \ref{directed-convergence-family} and the definition of $\mathcal{S}^*_X$-convergence, $\mathcal{S}^*_X$-convergence is equal to $\mathcal{S}^*_{\mathcal{D}X}$. Thus, $\mathcal{S}^*_X$-convergence  is topological iff $\mathcal{S}^*_{\mathcal{D}X}$ is topological iff  iff $\mathcal{D}X$ is a quasicontinuous space by Theorem \ref{T-qua}.
$\Box$

\vskip 3mm

\begin{rem} \rm 
Example \ref{example}  also shows that a $T_0$ space $X$ with a topological convergence class $\mathcal{S}^*_X$ may not be a quasicontinuous space. Since in $X = (\mathbb{N},co(\mathbb{N}))$, any subset $F =\ \ua F$, then $ F \leq G$ iff $G \subseteq F$. A directed family $\mathcal{F} $ of nonempty finite subsets must has a largest element $F$, i.e., $ \forall G \in \mathcal{F}, G \leq F \in \mathcal{F}$. Then $\mathcal{F} \to x$ iff $F = \{x\}$.  Thus, $((x_i)_I,x) \in \mathcal{S}^*_X$ iff $x$ is an eventual lower bound of $(x_i)_I$. Thus, $\mathcal{S}^*_X = \mathcal{C}_X$.
\end{rem}

\vskip 3mm

Let $X$ be a topological space and $A \subseteq X$. Denote $int_{\mathcal{D}X} A$ to be the interior of $A$ in space $\mathcal{D}X$.

\begin{prop}\rm  \label{inter}
Let $X$ be a finitely approximated space. Then
\begin{enumerate}
\item[(1)] For any $F \subseteq_f X$ and upper set $A$, $int_{\mathcal{D}X} \ua F =\ \Uparrow\! F$, $int_{\mathcal{D}X} A = \bigcup_{F \subseteq_f A}  \Uparrow\! F$
\item[(2)] $\{\Uparrow\! F: F \subseteq_f X\}$ form a base of the monotone determined topology $\mathcal{D}(X)$.
\end{enumerate}
\end{prop}

\proof 
(1) Suppose that $X$ is a finitely approximated space and $F \subseteq_f X$. Obviously, $\Uparrow\! F \subseteq \ua F$. Assume that $D$ is a directed subset of $X$ and $D \to y \in\ \Uparrow\! F$. By Proposition \ref{interpolation}, there exists some $G \subseteq_f X$ such that $F \ll_d G \ll_d y$. Then there exists some $d \in D$ such that $d \in \ua G$. Thus, $d \in\ \Uparrow\! F$ and then $D \cap \Uparrow\! F \not = \emptyset$. $\Uparrow\! F \in \mathcal{D}(X)$. 

Conversely, assume that $y \in int_{\mathcal{D}X} \ua F$ and $D$ is a directed set converging to $y$ with respect to $X$. Then $D$ is converges to $y$ with respect to $\mathcal{D}(X)$. It follows that there exists some $d \in D$ such that $d \in int_{\mathcal{D}X} \ua F$ and then $d \in \ua F$. Thus, $y \in\ \Uparrow F$. We conclude that $int_{\mathcal{D}X} \ua F =\ \Uparrow\! F$. 

We have $int_{\mathcal{D}X}A \supseteq \bigcup_{F \subseteq_f A} int_{\mathcal{D}X} \ua F = \bigcup_{F \subseteq_f A} \Uparrow F$. Conversely, assume that $y \in int_{\mathcal{D}X} A$. Since $X$ is a finitely approximated space, $fin(y)$ is a directed family and converges to $y$ with respect to $X$. By Lemma \ref{directed-convergence-family}, $fin(y) \to y$ with respect to $\mathcal{D}X$. Thus, there exists some $F \in fin(x)$ such that $F \ll_d y$ and $F \subseteq int_{\mathcal{D}X} A$. Therefore, $int_{\mathcal{D}X} A = \bigcup_{F \subseteq_f A} \Uparrow\! F$. 

(2) A direct conclusion from (1). $\Box$

\vskip 3mm

Obviously, a quasicontinuous space is a  finitely approximated space and for a monotone determined space $X$, $X$ is a quasicontinuous space iff it is a finitely approximated space. The following shows that a $T_0$ space is a finitely approximated space iff endowed with the monotone determined topology, it is a quasicontinuous space.

\vskip 3mm

\begin{thm}\rm  \label{qua-app}
Let $X$ be a $T_0$ space. Then $X$ is a finitely approximated space iff $\mathcal{D}X$ is a quasicontinuous space.
\end{thm}

\proof 
Since $F \ll_d x $ in $X$ iff $F \ll_d x$ in $\mathcal{D}X$ for any finite subset $F$ and $x \in X$, by Lemma \ref{directed-convergence-family}, $fin(x) = \{F :F \subseteq_f X, F \ll_d x \}$ is  a directed family that converges to $x$ with respect to $\mathcal{D}X$ iff it is a directed family converging to $x$ with respect to $X$. Therefore, $\mathcal{D}X$ is a quasicontinuous space iff $X$ is a finitely approximated space. $\Box$

\vskip 3mm

\begin{cor}\rm 
Let $X$ be a $T_0$ space. $\mathcal{S}^*_{X}$-convergence is topological iff $X$ is a finitely approximated space.
\end{cor}
\proof
By Theorem \ref{T-qua} and Theorem \ref{qua-app}. $\Box$

\vskip 3mm

Similarly with the case of quasicontinuous spaces, it is easy to verify that for a $T_0$ space $X$, $\mathcal{D}X$ is a continuous space iff for any $x \in X$, the following conditions are satisfied: (1) $\DDa_d x$ is a directed set; (2) $\DDa_d x$ converges to  $x$. Then by Theorem \ref{con-top}, given any $T_0$ space $X$, $\mathcal{S}_X$ is topological iff $X$ satisfying conditions (1) and (2).

In \cite[Theorem 4.17]{Z2022}, an equivalent condition for that $\mathcal{S}_X$ is topological is given as follows: $\mathcal{S}_X$-convergence is topological iff for all $x \in X$, the following hold: (1) $\DDa_d x$ is directed; (2) $\DDa_d x \to x$; (3) $\UUa_d x $ is open in $\mathcal{D}X$. Therefore, the condition (3) is redundant, which can be induced from (1) and (2).

\vskip 3mm

\begin{prop} \rm 
Let $X$ be a $T_0$ space satisfying that for any $x \in X$, $\DDa_d x$ is directed and $\DDa_d x \to x$. Then $\UUa_d x$ is open in $\mathcal{D}X$ for any $x \in X$.
\end{prop}
\proof
First, we show that if $x \ll_d y$, then there exists some $z \in X$ such that $x \ll_d z \ll_d y$. For any $y \in X$, consider $A = \{x \in X: \exists z \in X, x \ll_d z\ll_d y\}$. Assume that $x_1,x_2 \in A$ with $x_1 \ll_d z_1 \ll_d y$ and $x_2 \ll_d z_2 \ll_d y$. Then there exists some $z_3$ such that $x_1,x_2 \ll_d z_3 \ll_d y$ and then there exists some $x_3$ such that $x_1,x_2 \leq x_3 \ll_d z_3 \ll_d y$. Thus, $A$ is directed. We have $y \in cl(\DDa_d y) \in cl(A)$. Thus, $A \subseteq \DDa_d y$ is directed and converges to $y$. For any $x \ll_d y$, there exists some $x^\prime,z \in A$ such that $x \leq x^\prime \ll_d z \ll_d y$, that is, $x \ll_d z \ll_d y$.

Now, suppose that $D$ is a directed set converges to $y \in \UUa_d x$. Then there exists some $z \in \UUa_d x$ such that $z \ll_d y$ and there exists some $d \in D$ such that $z \leq d$. Thus, $d \in \UUa_d x$. It follows that $\UUa_d x$ is open in $\mathcal{D}X$.
$\Box$

\section{Lawson topology on monotone determined space}

In this section, we introduce the notion of quasi-liminf convergence and quasi-liminf topology. We will show that for a quasicontinuous space $X$, the quasi-liminf topology on $X$ is equal to the Lawson topology, i.e, the topology generated  by $\mathcal{O}(X) \vee \omega(X)$. Besides, we show that for a monotone determined space $X$, if the Lawson topology is compact, then $X$ is a dcpo endowed with the Scott topology.

\begin{defn} \rm
Let $X$ be a $T_0$ space. $x \in X$ and $(x_i)_I$ a net in $X$. $(x_i)_I$ is said quasi-liminf converging to $x$, denoted by $x \equiv_q \liminf x_i$, iff the following conditions hold:
\begin{enumerate}
\item[(1)] $((x_i)_I,x) \in \mathcal{S}^*_X$;
\item[(2)] for every eventual lower bound $y$ of $(x_i)_I, y \leq x$.
\end{enumerate}
\end{defn}

\vskip 3mm


\begin{lem} \rm \label{L* convergence}
Let $X$ be a $T_0$ space, $x \in X$ and $(x_i)_I$ a net in $X$. Then
the following statements are equivalent:
\begin{enumerate}
\item[(1)] $x \equiv_q \liminf y_j $ for all subnets  $(y_j)_J$ of $(x_i)_I$;

\item[(2)] $x \equiv_q \liminf x_i$ and $z \leq x $ if $z$ is a cofinal lower bound, i.e., if given any $i \in I$,
there exists $ j 
\geq i$ such that $z \leq  x_j$ .
\end{enumerate}
\end{lem}
 
\noindent{\bf Proof.}
(1) $\Rightarrow$ (2). Let $z$ be a cofinal lower bound of $(x_i)_I$.
Consider the subnet $(y_j)_J$ of $(x_i)_I$ consisting of all $x_i$ such that 
$ z \leq x_i $. Then $x$ is the quasi-liminf of $(y_j)_J$, 
 and $z$ is an eventual lower bound of this subnet. By definition of quasi-liminf convergence, $z \leq x$.

(2) $\Rightarrow$ (1). Let $(y_j )_J$ be a subnet of $(x_i)_I$ with $y_j = x_{f (j)}$. 
Let $z$ be an eventual lower bound of $(y_j)_J$. Then there exists $j_0$ such that $z \leq y_j $ for all $j_0 \leq j$. Then $\{f(j): j_0 \leq j \} $ is cofinal in $I$ and $z \leq x_{f(j)}$ for any $j_0 \leq j$. It follows that $z \leq x$ and then all eventual lower bounds of $(y_j)_J$ are contained in $\da x$.
Consider any quasi eventual lower bound $F$ of $(x_i)_I$. There exists $i_0 \in I$ such that $ x_i \in \ua F$ for all $ i\geq i_0$. There exists $j_0 \in J$ such that $i_0 \leq f(j)$ for all $j \geq j_0$. $F$ is a quasi eventual lower bound of $(y_j)_J$. Hence any directed family $\mathcal{F}$ of quasi eventual lower bounds of $(x_i)_I$ converges to $x$ is also a directed family of quasi eventual lower bounds of $(y_j)_J$.  Then $x \equiv_q \liminf y_j$. 
$\Box$

\vskip 3mm

\begin{defn}\rm 
Denote $\mathcal{L}^*$ the convergence class of all pairs of nets and elements in $X$ which satisfy the equivalent conditions of Lemma \ref{L* convergence}. We call the topology $\mathcal{L}^*(X) = \mathcal{L}^*_X(X)$  the quasi-liminf topology.
\end{defn}

\vskip 3mm

\begin{defn}  \rm 
Let $X$ be a monotone determined space. The Lawson topology on $X$, denoted by $\lambda(X)$, is the topology generated by $\mathcal{O}(X)$ and $\omega(X)$, where $\omega(X)$ is the lower topology on $X$ relative to the specialization order of $X$.
\end{defn}
\vskip 3mm

\begin{prop}\rm  \label{lambda quasilim}
Let $X$ be a $T_0$ space. Then the Lawson topology $\lambda(X)$ is coarser than $\mathcal{L^*}(X)$.
\end{prop}
\proof We need only to prove that any $U \in \mathcal{O}(X) \cup \omega(X)$ is open in $\mathcal{L}(X)$. 

(1) Since $\mathcal{L}_X^* \subseteq \mathcal{S}_X^*$,\ then $\mathcal{S}^*(X) \subseteq \mathcal{L}^*(X)$. By Lemma \ref{qs-d-topology}, $\mathcal{O}(X) \subseteq \mathcal{D}(X) \subseteq \mathcal{L}^*(X)$.

 
(2) 
Let $V =\ X \backslash\! \ua y$ for any $ y \in X$. Given any net $((x_i)_I,x) \in \mathcal{L}^*_X$ and $x \in V$, suppose that $(x_i)_I$ is not eventually in $V$. Then for any $i \in I$, there exists some $i \leq i^\prime$ such that $x_{i^\prime} \not \in V$. Consider the subnet of $(x_i)_I$ consisting of all $x_i$ such that $x_i \not \in V$, denoted by $(x_j)_J$. Then $((x_j)_J, x) \in \mathcal{L}^*_X$ by definition. Since $(x_j)_J \subseteq \ua y$, $y$ is a lower bound of $(x_j)_J$. Hence, $y \leq x$, a contradiction. Therefore, $(x_i)_I$ is eventually in $V$, i.e., $V$ is an open set relative to $\mathcal{L}^*(X)$.  

By (1) and (2), we conclude that $\lambda(X)$ is coarser than $\mathcal{L}^*(X)$.  $\Box$


\vskip 3mm

\begin{thm} \rm 
Let $X$ be a quasicontinuous space. Then $\mathcal{L}^*_X$-convergence is topological and agrees
with convergence in the Lawson topology.
\end{thm}
\proof
By the preceding proposition, $\lambda(X)$ is coarser than $\mathcal{L}^*(X)$. Thus, $\mathcal{L}^*_X$-convergence is contained in the convergence of $\lambda(X)$. We need only to show that for any net $(x_i)_I$ converging
to $x$ in the Lawson topology, $((x_i)_I,x) \in \mathcal{L}^*_X$. 
If $y \not \leq   x$, then $L\backslash\!\ua y$ is a Lawson open set containing $x$. There exists some $i_0 \in I$ such that $y \not \leq x_i$ for all $i \geq i_0$. It follows that all cofinal lower bounds of $(x_i)_I$ are lower than $x$. Since $X$ is quasicontinuous, $fin(x)$ is a directed family and converges to $x$. 
For any $F \in  fin(x)$, ${\Uparrow} F$ is open and hence Lawson open containing $x$. Since $(x_i)_I$ converges to $x$ with respect to the Lawson topology, then $(x_i)_I$ is eventually in ${\Uparrow} F$ and $F$ is an quasi eventual lower bound of $(x_i)_I$. Therefore, $((x_i)_I,x) \in \mathcal{L}^*_X$. $\Box$

\vskip 3mm

Finally, we investigate the compactness of the Lawson topology on a monotone determined space. We say that a monotone determined space is Lawson compact if the Lawson topology on it is compact and say that a subset is Lawson compact if it is compact relative to the Lawson topology. 

A space is called a d-space if it is a dcpo with respect to the specialization order and its topology is coarser than the Scott topology relative to the specialization order.

\vskip 3mm

\begin{lem}\rm \cite{ERNE2009} \label{dds}
Every monotone determined space $X$ which is also a d-space is a dcpo endowed with the Scott topology.
\end{lem}
\vskip 3mm

\begin{prop}\rm  \label{lawson dcpo}
Let $X$ be a monotone determined space. If $X$ is Lawson compact, then $X$ is a dcpo endowed with the Scott topology.
\end{prop}

\proof
Supposet that $X$ is Lawson compact. Given a directed subset $D$, $K=\bigcap_{d \in D}\ua d$ is a nonempty closed set since $\{\ua d\}_{d \in D}$ is a filtered family of closed sets. Hence $K$ is a non-empty Lawson compact set and then it is a compact set relative to the Scott topology. It follows that $K=\ \ua \text{min}~K$. For any $a,b\in \text{min}~K$, $\da a\ \cap \da b$ is Lawson compact, and $D$ is included in it.
It implies that $\bigcap_{d \in D}(\da a\ \cap \da b\ \cap\ua d)\neq\emptyset$. It is equivalent to $\da a\ \cap \da b\ \cap K\neq\emptyset$. This only happens when $a=b$. So $\text{min }\! K $ contains only one element. It means that  this element is the supremum of $D$. Therefore $X$ is a dcpo. Let $F$ be a closed set, then it is Lawson compact. For any directed set $D$ of $F$, we have  $\bigcap_{d\in D}(\ua d\cap F)\neq\emptyset$.
So $\bigvee D$ is included in  $F$.  It implies that $F$ is Scott closed. Then $F$ is a d-space. By Lemma \ref{dds}, $X$ is a dcpo endowed with the Scott topology. $\Box$

\vskip 3mm

\begin{prop}\rm \cite{Xi2016} \label{coherence patch}
For a $T_0$ space $X$, the following are equivalent.
\begin{enumerate}
\item[(1)] $X$ is well-filtered, compact and coherent. 
\item[(2)] $X$ is compact with respect to the patch topology, the topology with a closed subbase consisting of the
closed sets and the compact saturated sets.
\end{enumerate}
\end{prop}

\vskip 3mm

The following can be found in \cite[Lemma 4.1]{Xi2016} and \cite[Corollary 3.2]{Xi2016}.
\begin{lem}\rm \cite{Xi2016} \label{lemm}
Let $X$ be a dcpo endowed with the Scott topology. 
\begin{enumerate}
\item[(1)] If $X$ is well-filtered, then $X$ is coherent if $\ua x\ \cap \ua y$ is compact for all $x,y \in P$.
\item[(2)] If $X$ is Lawson compact then $X$ is  well-filtered.
\end{enumerate}
\end{lem}

\vskip 3mm

Then we get the following statement, which is an extension of \cite[Theorem 4.2]{Xi2016}.

\vskip 3mm

\begin{cor}\rm 
Let $X$ be a monotone determined space. Then the following statements are equivalent:
\begin{enumerate}
\item[(1)] $X$ is patch compact;
\item[(2)] $X$ is Lawson compact;
\item[(3)] $X$ is well-filtered, compact and ${\uparrow}x\cap{\uparrow}y$ is compact for any $x,y \in X$;
\item[(4)] $X$ is well-filtered, compact and coherent.
\end{enumerate}
\end{cor}





\section{Closures in locally hypercompact spaces}

Finally, we give a characterization of the closure of a subset in a quasicontinuous space.
Topological convergence classes and topologies on a set $X$ are one-to-one corresponding.\ When $\mathcal{E}$ is a topological convergence class,\ for any subset $F$ of $X$,\ its closure $\overline{F}= \{x \in X :(\xi,x) \in \mathcal{E},\ \xi \text{ is in } F\}$.\ This equality does not always hold for general convergence classes $\mathcal{E} \subseteq \{ (\xi,x): \xi \text{ is a net in } X,\ x \in X\}$.\ However,\ for some special $\mathcal{E}$,\ we can gain its closure by transfinite induction. 
\vskip 3mm

\begin{defn}\rm   \label{IND}
Let $P$ be a set.\ For any class $\mathcal{E} \subseteq \{ (\xi,x): \xi \text{ is a net in }P, x \in P\}$,\ any subset $F \subseteq P$,\ and any ordinal $\alpha \in ORD$,\ we define $F^{\alpha}$ and $F^{*}$ as follows:
\begin{align*}
F^{0}&=F  \\
F^{\alpha}&=\{ x \in P : \exists \xi \subseteq (\cup_{\beta < \alpha}F^{\beta}),\ (\xi,x) \in \mathcal{E} \} \\
F^{*}&= \cup_{\alpha \in ORD} F^{\alpha}
 \end{align*}

\end{defn}

\vskip 3mm

\begin{prop}\rm  \label{INDU}
Let $(P,\leq)$ be a poset and $\mathcal{E}$ be any class such that $\{ (\{y\},x): x,y \in P, \  x \leq y \} \subseteq \mathcal{E} \subseteq \{ (D,x): D \text{ is a directed subset of }P, \  x \in D^{\delta}\}$.\

\begin{enumerate}
\item[(1)] For any $F \subseteq P$,\ we have $\overline{F} = F^{*}$ in the topological space $(P,\mathcal{E}(P))$.

\item[(2)] The specialization order $\sqsubseteq$ of $(P,\mathcal{E}(P))$ is equal to $\leq$.

\item[(3)] $(P,\mathcal{E}(P))$ is a monotone determined space.
\end{enumerate}
\end{prop}

\proof
(1)\ It is easy to see that if $\alpha < \beta$, then $F^{\alpha} \subseteq F^{\beta}$,\ and if $F^{\alpha} = F^{\alpha + 1}$,\ then $F^{*} = F^{\alpha} $.\ For any subset $F$,\ there exists some $\alpha$ such that $F^{\alpha} = F^{\alpha + 1} = F^{*}$ since the cardinal of $F^{*}$ is less than or equal to that of $P$.\ We need only to show that $P \backslash F^{*}$ is open.\ 
For any $(D,x) \in \mathcal{E}$ such that $x \in P\backslash F^{*} $,\ if $D \subseteq F^{*}$,\ then $x \in F^{*}$,\ a contradiction.\ If $D \not\subseteq F^{*}$,\ then $D \cap (P \backslash F^{*}) \neq \emptyset$.\ By $\{(\{y\},x): x,y \in P, \  x \leq y \} \subseteq \mathcal{E}$,\ we know that $P \backslash F^{*} $ is an upper subset of $P$.\ Therefore,\ $D$ is eventually in $ P \backslash F^{*}$,\ i.e.,\ $P\backslash F^{*}$ is open.

(2)\ Given any $y \in P$,\ assume that $(D,x) \in \mathcal{E}$ and $x \in X \backslash\! \da y$.\ Then $D \cap (X \backslash\! \da y) \neq \emptyset$.\ Otherwise,\ $x \in D^{\delta} \subseteq\ \da y$,\ a contradiction.\ Therefore,\ $P \backslash\! \da y$ is open in $(P,\mathcal{E}(P))$ and then $x \sqsubseteq y$ implies $x \leq y$.\ If $x \leq y$, for any open subset $U$ such that $x \in U$,\ we have $y \in U$ by $ (\{y\},x) \in \mathcal{E}$.\ Thus,\ $x \leq y$ implies $x \sqsubseteq y$. 

(3)\ Let $U \in d(\mathcal{E}(P))$, i.e., for any directed subset $A$ such that $A \to x \in U$ relative to $\mathcal{E}(P)$, $A \cap U \not = \emptyset$. Then, for any $(D,x) \in \mathcal{E}$ such that $x \in U$, since $D \to x$ relative to $\mathcal{E}(P)$, $D \to x$ relative to $d(\mathcal{E}(P))$ by Theorem \ref{convergence}. It follows that $D \cap U \not = \emptyset$ and $D$ is eventually in $U$. Thus, $U$ is an open subset relative to $\mathcal{E}(P)$. $(P,\mathcal{E}(P))$ is a monotone determined space. $\Box$

\vskip 3mm

For a general convergence class $\mathcal{E}$, the conclusion of Proposition \ref{INDU} does not always hold.

\begin{exmp} \rm  
Let $\mathbb{N}^{\top}$ be the poset of natural numbers adding a top element $\top$,\ $\mathcal{E} = \{ (\mathbb{N},\top)\}$ and $A = \{2n + 1:\ n \in \mathbb{N}\}$. Then $A^{*} = A$.  Since $\mathbb{N} \to \top$ in $\mathcal{E}(\mathbb{N}^\top)$ and $A$ is a subnet of $\mathbb{N}$, $A \to \top$ relative $\mathcal{E}(\mathbb{N}^\top)$. Thus, the closure of $A$ must contain $\top$.
\end{exmp}
\vskip 3mm

From Proposition \ref{INDU}, we can also see that for any monotone determined space, the closure of its subsets can be gained by adding the limits of directed subsets of it by iteration. In particularly, for any quasicontinuous space $X$, the closure of any subset $A$ of $X$ can be gained by iterating the process only twice. Let $X$ be a monotone determined space. Considering the convergence class $DLim(X)$, then the topology of $X$ is determined by $DLim(X)$. By taking  $\mathcal{E}$ in Definition \ref{IND} as $DLim(X)$, we have the following statement.

\vskip 3mm

\begin{thm} \rm 
Let $X$ be a quasicontinuous space. Given any subset $A$ of $X$, $\overline{A} = A^2$.
\end{thm}
\proof Given any $x \in \overline{A}$, since $X$ is quasicontinuous, there exists a directed family
$$fin(x) = \{F \subseteq_f X: F \ll x \}$$
such that $\mathcal{F} \to x$. For any $F \in fin(x)$, since $x \in {\Uparrow} F \cap \overline{A}$ and ${\Uparrow} F = (\ua F)^\circ$ by Theorem \ref{d-quasicontinuity}, $(\ua F)^\circ \cap A \not = \emptyset$. Thus, $F\ \cap \da A \not = \emptyset$. 

Now, we verify that $\{F\ \cap \da\! A: F \in fin(x)\}$ is a directed family. Since $fin(x)$ is directed, for any $F_1,F_2 \in fin(x)$, there exists an $F \in fin(x)$ such that $F \subseteq\ \ua F_1\ \cap \ua F_2$. Then, for any $x \in F\ \cap \da A$, there exists an $a \in A$ and $k_i \in F_i$ such that $k_i \leq x \leq a$  for $i =1,2$. Therefore, $F\ \cap \da A \subseteq\ \ua (F_1\ \cap \da A)\ \cap\ \ua (F_2\ \cap \da A)$, i.e., $\{F\ \cap \da A: F \in fin(x)\}$ is a directed family. 
There exists a directed subset $D \subseteq \bigcup_{F \in fin(x)}(F\ \cap \da A)$ such that $D \cap (F\ \cap \da A) \not= \emptyset$ holds for any $F \in fin(x)$  by Ruddin's Lemma. Since $X$ is quasicontinuous, for any open subsets $U$ with $x \in U$, there exsits some $F \in fin(x)$ such that $F \subseteq U$. Thus, $D \cap U \not = \emptyset$. Pick $x_U \in U \cap D \subseteq\ \da A \subseteq A^1$ for each open neighborhood $U$ of $x$. Then, all $x_U$ form a net that converges to $x$. Thus, $x \in A^2$. It follows that $\overline{A} = A^2$.

\end{document}